\documentclass[reqno,10pt,a4paper]{amsart}
\usepackage{graphicx,color,amssymb,latexsym,amsfonts,amsmath}
\usepackage{mathrsfs}
\usepackage{graphics}
\usepackage{lscape}
\newtheorem{lemma}{\bf Lemma}[section]

\newtheorem{corollary}{\bf Corollary}
\newtheorem{remark}{Remark}

\newtheorem{theorem}{\bf Theorem}[section]

\numberwithin{equation}{section} \theoremstyle{plain}
\theoremstyle{definition}

\DeclareMathOperator*{\argmax}{argmax}

\input cyracc.def

\begin{document}
\title[Likely path to extinction
in simple branching models] {Likely path to extinction   in simple
branching models with large initial population}
\author{Klebaner F.C.}
\address{School of Mathematical Sciences, Building 28M, Monash
University, Clayton Campus, Victoria 3800, Australia.}
\email{fima.klebaner@sci.monash.edu.au}

\author{Liptser R.}
\address{School of Mathematical Sciences, Building 28M, Monash
University, Clayton Campus, Victoria 3800, Australia.}
\email{rliptser@sci.monash.edu.au, liptser@eng.tau.ac.il}
\thanks{Research supported by the Australian Research Council Grant
DP0451657} 
 \maketitle

\begin{abstract}
We give explicit formulae for   most likely paths to extinction in
simple branching models when
  initial population is large. In
discrete time we study the
  Galton-Watson process and in continuous time the Branching
diffusion.
  The most likely paths are found with the help of the Large
  Deviation Principle
(LDP). We also find asymptotics for the extinction probability,
 which gives a new expression in continuous time and recovers
 the known formula in discrete time. Due to the
non-negativity of the processes, the proof of LDP at the point of
extinction uses a nonstandard argument of independent interest.
\end{abstract}

\section{\bf Introduction and main results}
\label{sec-1}

In population genetics it is often important to look back at the
development of   populations. In this paper we consider the
question of how extinctions occur, and in particular, what path a
population takes on the road to extinction. Using asymptotic
analysis when initial population values are large,  we are able to
find most likely path to extinction as well as the extinction
probability in two simple
  branching models in discrete  and continuous time.
 In both  examples  we use the large deviation principle (LDP)
which is non-standard
since random processes are nonnegative, and we use
trajectories ending up at zero.

One of the contributions of this paper is in rigorous proofs of
the LDP for processes on half space. It may appear to the reader
that the LDP follows from known results in Markov chains and
diffusions. This is only partly correct.  The standard proof of
the lower bound in the local LDP relies on the change of measure.
This requires a certain point (the point where maximum in the
Fenchel-Legendre transform is achieved) to be finite. In our case
this point is at infinity, breaking down the standard approach. We
therefore give complete proofs of LDP's  in Sections \ref{LDPdisc}
(discrete time) and
 \ref{LDPcont} (continuous time) following the
 scheme of Puhalskii \cite{puh2}.
 His approach states that the LDP is equivalent to exponential tightness plus local
 LDP,  and is  based on the stochastic exponential method
  (rather than the
Laplace transform). Although we follow the   scheme of Puhalskii
\cite{puh2} we do not use idempotent probability and give direct
proofs.
 Since these proofs are more technical, we placed them at the end, after results on
 extinction.
 Once the LDP is established, the problem of finding most likely path to extinction
 is in effect the problem of minimization of the rate function.
 This is typically a difficult problem due to nonlinearity.
  We are able to solve it by setting up  the Bellman equation in discrete case,
 Section \ref{sec-2}, and
  a dynamical control problem in continuous case,
Section \ref{sec-3}.

\subsection{Galton-Watson process.}

A prototype of a branching model in discrete time is the
Galton-Watson process,   described as follows.

  Let $X_n$ denote the population size at  time
$n$, and $\xi^j_{n+1}$ the number of offspring of the $j$th
individual. For each $n=1,2,\ldots$, $\{(\xi^j_n)_{j\ge 1}\} $ is
the sequence of independent identically distributed integer-valued
random variables with the probability distribution function
$
\mathsf{P}(\xi_n^j=\ell)=p_\ell, \ \ell=0,1,\ldots.
$
The population size at time $n+1$ is given by
$$
X_{n+1}=\sum_{j=1}^{X_n}\xi^j_{n+1},
$$
where  $X_0=K>0$. The state $\{0\}$ is absorbing, and the
branching process $(X_n)_{n\ge 0}$ might be absorbed in $\{0\}$ at
the extinction time
$$
\tau=\inf\big\{n:X_n=0\big\}.
$$
If $p_0=0$, the population does not become extinct. However if
$p_0>0$, it is well known (e.g.  \cite{Har}, \cite{AthN}) that the
extinction time $\tau$ is finite with probability one if and only
if the offspring mean $ \frak{m}=\sum_{\ell\ge 1}\ell p_\ell $
does not exceed  one ($\frak{m}\le 1$). Moreover, for any
$\frak{m}$, the distribution function of $\tau$ is computed using
the offspring probability generating function $
\frak{f}(s)=\sum_{\ell\ge 0} p_\ell s^\ell,$  $0\le s \le 1$: for
any $N\ge 1$,
\begin{equation}\label{1.2ba}
\mathsf{P}(\tau\le N)=(\frak{f}_N(0))^K,
\end{equation}
where $\frak{f}_n(s)$ is the $n$-th iterate of $\frak{f}(s)$, i.e.
$\frak{f}_n(s)=\frak{f}(\frak{f}_{n-1}(s))$ with
$\frak{f}_1(0)=\frak{f}(0)=p_0$.

A  natural  question is how to find  the ``path to extinction''
  given that extinction occurred at time $N$, $\tau=N.$
  The conditional distribution of the chain conditioned on extinction: for $n=1,\ldots,N-1$,
$$
\pi_{n|N}(i):=\mathsf{P}(X_n=i|\tau=N), \ i=1,2,\ldots
$$
gives the complete description. It can be used to find the
conditional median or the traditional optimal  in the mean square
sense
 estimate $\widehat{X}_n=\sum_{i=1}^\infty
i\pi_{n|N}(i)$. Unfortunately such computations are   involved,
even using the Markov property of $(X_n)$.
However, for large
values of $X_0=K$, one path has an overwhelmingly large
probability compared to the rest. Consider the normed branching
process $$ x^K_n=\frac{X_n}{K}.$$ The limit in probability
$\mathsf{P}\text{-}\lim_{K\to\infty}x^K_n=\hat{x}_n$ exists (see
\cite{Kleban},  \cite{KlNer}) and satisfies $
\hat{x}_{n+1}=\frak{m}\hat{x}_n, \quad \hat{x}_0=1.$ The process
$\hat{x}_n$ is always positive, irrespective of the value of
$\frak{m}$, so that, the approximation $\hat{x}_n$ is inadequate for study
of extinction, the fact is
already mentioned in \cite{Barbour}.
 In the approach we take, $(x^K_n)_{n\le N}$ is
approximated on the set $\{\tau\le N\}$ by a deterministic
sequence $u^*_\cdot:=(u^*_n)_{n\le N}$ with $u^*_0=1$, positive
$u^*_n$'s and $u^*_N=0$,  such that for small $\delta>0$ and large
$K$,
\begin{equation*}
\mathsf{P}\Big(\sum_{n=1}^N|x^K_n-u^*_n|\le \delta\Big)\approx
\mathsf{P}\Big(\tau\le N\Big) .
\end{equation*}
This choice of $u^*_\cdot$ might be warranted by the following
argument. Since $\frak{f}_n(0)$ increases in $n$, for
 large $K$,  $(\frak{f}_N(0))^K$ is considerably
larger than any of $(\frak{f}_n(0))^K$ for $n<N$. Then, by
\eqref{1.2ba},
$
\mathsf{P}(\tau\le N)=\mathsf{P}(\tau= N)+\mathsf{P}(\tau\le
N-1)\approx \mathsf{P}(\tau= N).
$
Consequently, for any $u_\cdot=(u_n)_{n\le N}$ with $u_0=1$ and
$u_n\ge 0$,
\begin{equation*}
\mathsf{P}\Big(\sum_{n=1}^N|x^K_n-u_n|\le \delta\Big)\lesssim
\mathsf{P}\Big(\tau\le N\Big) .
\end{equation*}
For large $K$, extinction for the process $x_n^K$ is a rare event,
since the limit process $\hat{x}_n$ is positive. Therefore, as in
\cite{KL}, we approach the problem of extinction using the large
deviations theory, obtaining a new result  as well as recover an
asymptotic version of the well-known result \eqref{1.2ba} by using
this theory. According to LDP, Theorem \ref{theo-A.1} and by
analogy with the maximal likelihood estimator, the path
$(u^*_n)_{n\le N}$ is said to be the most likely path to
extinction of the normed population $x^K_n$.

Clearly, $\tau$ is
the extinction time  for both processes $X_n$ and $x^K_n$, so
that, $Ku^*_n$ (with large $K$) sets the pattern for the
extinction path in the original branching process.

Figure below demonstrates likely paths to   extinction   for a
binary splitting model with different parameters, $p=p_0$, illustrating the
general result.
\begin{figure}[ht]
\begin{center}
\includegraphics[angle=0,width=4.5in,height=2.5in]{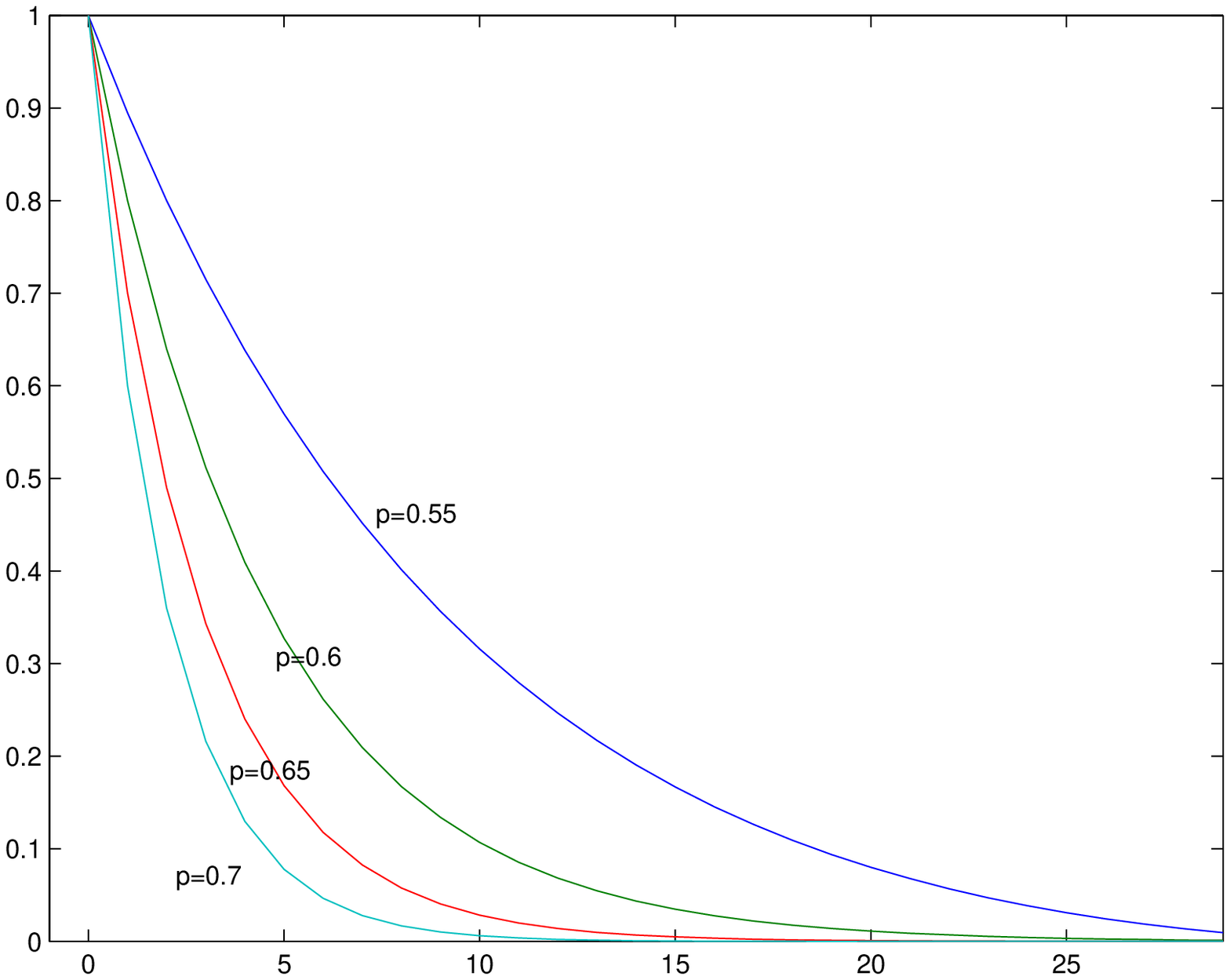}
\end{center}
\end{figure}

For formulating the main result,  we use the log moment generating
function, assuming its existence up to some $t_0>0$,
\begin{equation}\label{log}
\frak{g}(t)=\log \sum_{\ell\ge 0}e^{t\ell}p_\ell, \
t\in(-\infty,t_0).
\end{equation}
It is related to the moment generating function by
$$
\log\frak{f}_n(0)\equiv\frak{g}_n(-\infty) \quad\text{(Lemma
\ref{prop-3.1a})}.
$$
\begin{theorem}\label{theo-main1}
Assume $p_0>0$ and \eqref{log}. Then, for any $N\ge 1$,

{\rm (i)}
$$
\begin{aligned}
(u^*_n)_{n\le N}&=\argmax_{\substack{u_0=1,u_N=0\\
u_n>0,n\le N-1}}\lim_{\delta\to 0}\lim_{K\to\infty}\frac{1}{K}\log
\mathsf{P}\Big(\sum_{n=1}^N|x^K_n-u_n|\le \delta\Big)
\end{aligned}
$$
 with
\begin{equation}\label{MLPDisc}
u^*_n=\prod_{1\le i\le n}\frak{g}'(\frak{g}_{N-i}(-\infty)), \
n\le N,
\end{equation}
where $\frak{g}_i(t)$ is $i$-th iterate of $\frak{g}(t)$,
$\frak{g}_0(t)=t$.

{\rm (ii)}
$$
\lim_{\delta\to 0}\lim_{K\to\infty}\frac{1}{K}\log
\mathsf{P}\Big(\sum_{n=1}^N|x^K_n-u^*_n|\le
\delta\Big)=\lim\limits_{K\to\infty}
\frac{1}{K}\log\mathsf{P}(\tau\le n)\big).
$$
\end{theorem}

\subsection{Branching diffusion.}

  In  continuous time, we
consider the model of a branching diffusion $X_t$ defined by the
It\^o equation
\begin{equation}\label{KKK}
dX_t=\alpha X_tdt+\sigma\sqrt{X_t}dB_t
\end{equation}
with a positive initial condition $X_0=K$, where $B_t$ is a
Brownian motion,    $\sigma^2>0$, and $\alpha\in\mathbb{R}$.
Stochastic equation \eqref{KKK} possesses  a strong nonnegative
solution. Since the diffusion parameter degenerates, one way to
see this is to construct the solution from the following
approximating sequence $(X^i_t)_{ i\ge 1}$:
$$
X_t:=X^1_tI_{\{t\le \tau_1\}}+\sum_{i\ge
1}X^i_{\tau_i}I_{\{\tau_i<t\le\tau_{i+1}\}},
$$
where $ dX^i_t=\alpha X^i_tdt+\sigma\sqrt{|X^i_t|\vee i^{-1}}dB_t,
$ $X^i_0=K$, and $\tau_i=\inf\{X^i_t\le i^{-1}\}$ the increasing
sequence of stopping times $(\tau_i)_{i\ge 1}$ relative to the
filtration generated by Brownian motion $(B_t)$ (see also Theorem
13.1, \cite{Kleb}). The strong uniqueness of \eqref{KKK} follows
from Yamada-Watanabe's theorem (see, e.g. Rogers and Williams, p.
265 \cite{RW}) since its drift and diffusion parameters are
Lipschitz and H\"older (with coefficient $\frac{1}{2}$) continuous
respectively.

Obviously,
$$
\tau=\inf\{t:X_t=0\}=\lim_{i\to\infty}\tau_i.
$$
We  analyze  the normed process $x^K_t=\frac{X_t}{K}$. Due to
\eqref{KKK},  $x^K_t$ solves the It\^o equation
\begin{equation}\label{nd}
dx^K_t=\alpha x^K_tdt+\frac{\sigma}{\sqrt{K}}\sqrt{x^K_t}dB_t,
\end{equation}
with $x^K_0=1$. It can be readily shown that
$\mathsf{P}\text{-}\lim_{K\to\infty}x^K_t=\hat{x}_t$ exists and
solves $ \dfrac{d\hat{x}_t}{dt}=\alpha\hat{x}_t, \ \hat{x}_0=1. $
However, $\hat{x}_t$ is always positive and is far from to be
estimated path to extinction. As in the discrete time, in order to
evaluate path to extinction for $(x^K_t)_{t\le T}$ for fixed
$T>0$, we approximate $(x^K_t)_{t\le T}$ on the set $\{\tau\le
T\}$ by a deterministic function $(u^*_t)_{t\le T}$ with $u^*_0=1,
u^*_T=0$ and $u^*_t>0$, such that for a small $\delta>0$ and large
$K$,
\begin{equation*}
\mathsf{P}\Big(\sup_{t\le T}|x^K_t-u^*_t|\le \delta\Big)\ge
\mathsf{P}\Big(\sup_{t\le T}|x^K_t-u_t|\le \delta\Big)
\end{equation*}
for any $(u_t)_{t\le T}$ from the set $\{u_0=1, (u_t>0)_{t< T},
u_T=0\}$.

Unfortunately, the useful formula of  \eqref{1.2ba} type  is not
known to us in this case. Here we obtain its asymptotic version as
$K\to\infty$,  see (ii) below.

\begin{theorem}\label{theo-main2} For any $T>0$,
\begin{enumerate}
\item[(i)] $$
(u^*_t)_{t\le T}=\argmax_{\substack{u_0=1,u_T=0\\
u_t>0,t<T}}\lim_{\delta\to 0}\lim_{K\to\infty}\frac{1}{K}\log
\mathsf{P}\Big(\sup_{t\le T}|x^K_t-u_t|\le \delta\Big)
$$
 is given by
\begin{equation}\label{MLPcont}
u^*_t=\begin{cases}
  e^{-\alpha t}\Big(1-\frac{1-e^{-\alpha t }}{1-e^{-\alpha T }}\Big)^2    , &
 \alpha\ne 0 \\
     \big(1-\frac{t}{T}\big)^2 , & \alpha=0.
   \end{cases}
\end{equation}
\item[(ii)]
\begin{gather*}
\lim\limits_{K\to\infty}
\frac{1}{K}\log\mathsf{P}(\tau\le T)
=\lim_{\delta\to 0}\lim_{K\to\infty}\frac{1}{K}\log
\mathsf{P}\Big(\sup_{t\le T}|x^K_t-u^*_t|\le \delta\Big)
\\
=-\begin{cases}
\frac{1}{\sigma^2}\frac{\alpha}{1-e^{-\alpha T}}, & \alpha\ne 0
\\
\frac{1}{\sigma^2T}, & \alpha=0
\end{cases}
\end{gather*}
\end{enumerate}
\end{theorem}

\begin{corollary}
\begin{enumerate}
\item $u^*_\cdot$ has the remarkable  property{\rm :} it is the same
for subcritical and supercritical case: $u^*_t(\alpha)\equiv
u^*_t(-\alpha)$. \item For  large $K$, the probability of
extinction in $[0,T]$ is given by
\begin{equation*}
\mathsf{P}(\tau\le
T)\approx\exp\Big(-\frac{K}{\sigma^2}\frac{\alpha}{1-e^{-\alpha
T}} \Big).
\end{equation*}
In particular, for $\alpha=0$, $ \mathsf{P}(\tau\le T)\approx
e^{-\frac{K}{2\sigma^2T}}. $
\end{enumerate}
\end{corollary}

\section{\bf Proof of Theorem \ref{theo-main1}}
\label{sec-2}

We begin with
\begin{lemma}\label{prop-3.1a}
For any $n\ge 1$, $ \frak{g}_n(-\infty)=\log \frak{f}_n(0). $
\end{lemma}
\begin{proof}
The result follows by   induction from the identity
$
\frak{g}_n(\log t)\equiv\log \frak{f}_n(t)
$
for $t\in (0,t_0).$
Write
$$
\frak{g}(\log(t))=\log\sum_{\ell\ge 0}e^{\ell\log(t)}p_\ell =\log
\sum_{\ell\ge 0}e^{\log(t^\ell)}p_\ell =\log \sum_{\ell\ge
0}t^\ell p_\ell=\log\frak{f}(t).
$$
If $\frak{g}_{n-1}(\log t)\equiv\log \frak{f}_{n-1}(t)$, then
$$
\frak{g}_n(\log
t)=\frak{g}(\frak{g}_{n-1}(t)=\frak{g}(\log(\frak{f}_{n-1}(t))=
\log\frak{f}(\frak{f}_{n-1}(t))=\log(\frak{f}_n(t)).
$$
\end{proof}

The proof of Theorem \ref{theo-main1} is done in a number of
steps.

\medskip
\textbf{(1)} Recall that $\frak{g}(t)$ is convex function with
$\frak{g}(0)=0$, $\frak{g}(-\infty)=\log(p_0)$  and
$\frak{g}'(t)>0, \ t>-\infty$ while
$\frak{g}'(-\infty)=\lim_{t\to\infty}\frak{g}'(t)=0$.

\medskip
\textbf{(2)} By the local LDP (see, Theorem \ref{theo-A.1}), for $
u_0=1, u_N= 0 $ and other positive $u_n$'s, it holds
\begin{equation*}
\lim_{\delta\to
0}\lim_{K\to\infty}\frac{1}{K}\log\mathsf{P}\Big(\sum_{n\le N}
|x^K_n-u_n|\le\delta\Big)=-\sum_{n\le N}I(u_n,u_{n-1}).
\end{equation*}

\smallskip
\textbf{(3)} In order to find $(u^*_n)_{n\le N}$ such that for
$u_0=1,u_n>0, u_N=0$
\begin{equation}\label{3.1de}
\sum_{i\le n}I(u_i,u_{i-1})\ge \sum_{i\le n}I(u^*_i,u^*_{i-1}),
\end{equation}
we apply the Dynamic Programming.

Since $u_N=0$,
\begin{equation}\label{3.3bc}
I(u_N,u_{N-1})=\sup_{t\in(-\infty,t_0)}(-u_{N-1}\frak{g}(t))=-u_{N-1}\frak{g}(-\infty)
=:B_n(u_{N-1})
\end{equation}
is the boundary condition for the Bellman equation
\begin{equation}\label{BelEq}
B_{n}(u_{n-1})=\inf_{u>0}\Big[B_{n+1}(u)+I(u,u_{n-1})\Big], \ 1\le
n\le N-1.
\end{equation}
For $n=N-1$, we have
\begin{equation}\label{3.2aa}
B_{N-1}(u_{N-2})=\inf_{u>0}\Big[-u\frak{g}(-\infty)+\sup_{t\in(-\infty,t_0)}
\{tu-u_{N-2}\frak{g}(t)\}\Big].
\end{equation}
\eqref{3.2aa} provides the inequality,
$$
B_{N-1}(u_{N-2})\ge \inf_{u>0}\big[-u\frak{g}(-\infty)+
tu-u_{N-2}\frak{g}(t)\big], \  \forall \ t\in (-\infty,t_0)
$$
which, with  $t=\frak{g}(-\infty)$, is transformed into
\begin{equation}\label{BN1}
B_{N-1}(u_{N-2})\ge -u_{N-2}\frak{g}_2(-\infty).
\end{equation}
We show that the above inequality is  equality. For $u,u_{N-2}>0$,
``$\sup_t$'' in \eqref{3.2aa} is attained at the point
$t^*=t^*(u,u_{N-2})$, so that, for any $u>0$,
\begin{equation*}
B_{N-1}(u_{N-2})\le u\big[t^*(u,u_{N-2})-\frak{g}(-\infty)\big]-
u_{N-2}\frak{g}(t^*(u,u_{N-2})).
\end{equation*}

We choose $u=u^*_{N-1}$ such that
$t^*(u^*_{N-1},u_{N-2})=\frak{g}(-\infty)$. This is possible since
$$
\begin{aligned}
&
\frak{g}(-\infty)=\log p_0,\quad t^*(0,u_{N-2})=-\infty,
\\
&
\frak{g}'(-\infty)=0,\quad t^*(\frak{m},u_{n-2})=0,\quad \frak{g}'(0)=\frak{m},
\end{aligned}
$$
so that, the existence of $u^*_{N-1}$
follows from continuity, in $u$, of $t^*(u,u_{N-2})$.

The choice of $u^*_{N-1}$ gives the inequality
$$
B_{N-1}(u_{N-2})\le -
u_{N-2}\frak{g}(t^*(u^*_{N-1},u_{N-2}))=\frak{g}\big(\frak{g}(-\infty)\big)=
\frak{g}_2(-\infty).
$$
Consequently, the opposite inequality for \eqref{BN1} holds true
and, therefore,
$$
B_{N-1}(u_{N-2})=-u_{N-2}\frak{g}_2(-\infty).
$$
It is obvious too that for any $u_{N-2}>0$,
$$
u^*_{N-1}=u_{N-2}\frak{g}'(t^*(u^*_{N-1},u_{N-2}))=u_{N-2}\frak{g}'(\frak{g}(-
\infty)).
$$

Further, by induction, we find the following pairs:
\begin{equation*}
\begin{aligned}
u^*_{N-1}&=\frak{g}'\big(\frak{g}(-\infty)\big)u^*_{N-2}
\\
B_{N-1}(u^*_{N-2})&=-\frak{g}_2(-\infty)u^*_{N-2}
\\
&\ldots
\\
u^*_{N-2}&=\frak{g}'\big(\frak{g}_2(-\infty)\big)u^*_{N-3}
\\
B_{N-2}(u^*_{N-3})&=-\frak{g}_3(-\infty)u^*_{N-3}
\\
&\vdots
\\
u^*_1&=\frak{g}'\big(\frak{g}_{n-1}(-\infty)\big)u_0
\\
B_1(u_0)&=-\frak{g}_n(-\infty)u_0 \ (u_0=1).
\end{aligned}
\end{equation*}
With chosen $(u^*_n)_{1\le n\le N-1}$, the Bellman equation
\eqref{BelEq} is transformed into the backward recurrent equation
\begin{equation*}
B_{n}(u^*_{n-1})=B_{n+1}(u^*_n)+I(u^*_n,u^*_{n-1}), \ 1\le n\le
N-1
\end{equation*}
with boundary condition $-u^*_{N-1}\frak{g}(-\infty)$ (see,
\eqref{3.3bc}).

\smallskip
Thus, $ B_1(1)=\sum_{1\le n \le N}I(u^*_n,u^*_{n-1}). $

On the other hand, the Bellman equation also yields
$$
B_1(1)\ge \sum_{1\le n \le N-1}I(u_n,u_{n-1})+B_N(u_{N-1})=
\sum_{1\le n \le N}I(u_n,u_{n-1})
$$
what proves \eqref{3.1de}.

\medskip
\textbf{(4)} We recall that $\sum_{n\le N}I(u^*_n,u^*_{n-1})
=-\frak{g}_n(-\infty)$, that is, by Lemma \ref{prop-3.1a} and
\eqref{1.2ba},
$$
\sum_{1\le n\le
N}I(u^*_n,u^*_{n-1})=-\log\frak{f}_N(0)=-\frac{1}{K}\log\mathsf{P}(\tau\le
N), \ \forall \ K>0.
$$

\medskip
\textbf{(5)} Thus, \textbf{(1)}-\textbf{(3)} imply the statement
(i); formula \eqref{MLPDisc} follows from   recurrence
 $
u^*_{n}=\frak{g}'\big(\frak{g}_2(-\infty)\big)u^*_{n-1} $,
$u^*_0=1$.

Finally (ii) follows from \textbf{(4)}. \qed

\section{\bf Proof of Theorem \ref{theo-main2}}
\label{sec-3}

We apply the LDP Theorem \ref{theo-B.1}. By the local LDP, with $u_0=1$,
$u_t>0$ and $u_T=0$, we have
$$
\lim_{ \delta\to
0}\lim_{K\to\infty}\frac{1}{K}\log\mathsf{P}\Big(\sup_{t\le T}
\big|x^K_t-u_t\big|\le\delta\Big)=-J_T(u),
$$
where $J_T(u)=
   \begin{cases}
     \frac{1}{2\sigma^2}\int_0^T\frac{(\dot{u}_t-u_t)^2}{u_t}I_{\{u_t>0\}}dt , &
 u_0=1,du_t=\dot{u}_tdt
      \\
     \infty , & \text{otherwise}.
   \end{cases}
$

Therefore (i) is reduced to  minimization of $J_T(u)$ in a class
of absolutely continuous test functions $u_t$ with $u_0=1$, $u_t>0$ and
$u_T=0$.

Set $w_t=\dfrac{\dot{u}_t-u_t}{\sqrt{u_t}}$, $t\in[0,T)$ and notice that
the minimization of $J_T(u_\cdot)$ is equivalent to the following control problem
with the controlled process $u_t$, solving a differential equation
$$
\dot{u}_t=\alpha u_t+\sqrt{u_t}w_t, \quad t\in [0,T)
$$
subject to $u_0=1$.   The control action $w_t$ belongs to a class
of measurable functions with $\int_0^Tw^2_tdt<\infty$ bringing
$u_t$ to zero at the time $T$. The control action $w^*_t$ from
this class is optimal if for any $w_t$,
\begin{equation*}\
\int_0^T(w^*_t)^2dt\le \int_0^Tw^2_tdt.
\end{equation*}
If $w^*_t$ exists, then the controlled process $u^*_t$ related to $w^*_t$ minimizes
$J_T(u_\cdot)$ in the required class of continuous functions $u_\cdot=(u_t)_{t\le T}$.

In order to find $w^*_t$, it is convenient to deal with (recall
$u_t\ge 0$) $v_t=\sqrt{u_t}$ since $v_t$ solves the linear
differential equation $
\dot{v}_t=\frac{\alpha}{2}v_t+\frac{1}{2}w_t,\quad v_0=1. $ If
$w^*_t$ exists, then $w^*_t$ brings $v_t$ to zero at the time $T$,
that is, $
0=v_T=e^{\frac{\alpha}{2}T}+\int_0^Te^{\frac{\alpha}{2}(T-t)}w^*_tdt
$ or, equivalently,
\begin{equation}\label{bla}
-1=\frac{1}{2}\int_0^Te^{-\frac{\alpha}{2}t}w^*_tdt.
\end{equation}
Hence, by the Cauchy-Schwarz inequality $ 1\le\frac{1}{2}
\int_0^Te^{-t\alpha}dt\int_0^T(w^*_t)^2dt$, that is, the following
lower bound holds: $
\int_0^T(w^*_t)^2dt\ge\dfrac{2\alpha}{1-e^{-\alpha T}}. $ This
lower bound is valid for any $w_t$ providing \eqref{bla} , so
that, the condition
$$
\int_0^T(w^*_t)^2dt=\dfrac{2\alpha}{1-e^{-\alpha T}}
$$
is valid
for
$w^*_t=ce^{-t\frac{\alpha}{2}}$ for any constant $c$, bring
$w^*_t=c^*e^{-t\frac{\alpha}{2}}$ with $c^*$ solving
$$
-1=\int_0^Te^{-t\frac{\alpha}{2}}w^*_tdt=c^*\int_0^Te^{-t\alpha}dt.
$$
Hence,
\begin{gather*}
c^* =
   \begin{cases}
     -\frac{2\alpha}{1-e^{-T\alpha}}, & \alpha\ne 0
     \\
     -\frac{2}{T} , & \alpha=0
   \end{cases}
   \quad\text{and}\quad
   w^*_t=
   \begin{cases}
-\frac{2\alpha e^{-t\frac{\alpha}{2}}}{1-e^{-T\alpha}} , &
\alpha\ne 0
\\
-\frac{2}{T}  , & \alpha=0
   \end{cases}
   \\
\int_0^T(w^*_t)^2dt=
   \begin{cases}
     \frac{2\alpha}{1-e^{-\alpha T}} , & \alpha\ne 0 \\
     \frac{2}{T} , & \alpha=0.
   \end{cases}
\end{gather*}
Finally, we find that
\begin{equation*}
\begin{aligned}
v^*_t&=e^{t\frac{\alpha}{2}}-\frac{\alpha}{1-e^{-T\alpha}}
\int_0^te^{(t-s)\frac{\alpha}{2}}e^{-s\frac{\alpha}{2}}ds
\\
&=e^{t\frac{\alpha}{2}}\Big[1-\frac{1-e^{-t\alpha}}{1-e^{-T\alpha}}\Big]
=e^{t\frac{\alpha}{2}}\Big(\frac{e^{-t\alpha}-e^{-T\alpha}}{1-e^{-T\alpha}}\Big)
\end{aligned}
\end{equation*}
and, since $u^*_t=(v^*_t)^2$, we obtain \eqref{MLPcont}
and the proof of (i) is complete.

\medskip
(ii) By (i),
\begin{equation}\label{TTT}
J_T(u^*)=\frac{1}{\sigma^2}\frac{\alpha}{1- e^{-\alpha T}}
\end{equation}
We show  that
\begin{equation*}
\lim_{K\to\infty}\frac{1}{K}\log\mathsf{P}\big(\tau\le
T\big)=J_T(u^*).
\end{equation*}
To this end, use the fact that $\{\tau\le T\}=\{(\omega,t):\exists
t\le T, \ x^K_t(\omega)=0\}$. For notational convenience denote
$\frak{A}:=\{\tau\le T\}$. Set $\frak{A}^{\sf cl}$ and
$\frak{A}^{\sf int}$ the  closure and interior of $\frak{A}$.
Then, by the LDP, we have
\begin{equation*}
\begin{aligned}
\varlimsup_{K\to\infty}\frac{1}{K}\log\mathsf{P}\big(\frak{A}^{\sf
cl}\big)
\le -\inf_{u:\left\{\substack{u_s>0,s<t;\\
u_t=0\\
t\le T}\right.}J_t(u)=-\inf_{t\le T}J_t(u^*)
\\
\varliminf_{K\to\infty}\frac{1}{K}\log\mathsf{P}\big(\frak{A}^{\sf
int}\big)
\ge -\inf_{u:\left\{\substack{u_s>0,s<t;\\
u_t=0\\
t\le T}\right.}J_t(u)=-\inf_{t\le T}J_t(u^*).
\end{aligned}
\end{equation*}
Since $\varliminf_{K\to\infty}=\varlimsup_{K\to\infty}$ implies
the existence of $\lim_{K\to\infty}$, it remains to show that
$\inf_{t\le T}J_t(u^*)=J_T(u^*)$.

Notice that \eqref{TTT} is valid with $T$ replaced by any $t<T$
with $u^*_\cdot$ replaced by the corresponding  $u^{*,t}_\cdot=
\{u^{*,t}_0=1;u^{*,t}_s>0,s<t;u^{s,t}_t=0\}$. In other words, for
any $t$,
$$
J_t(u^{*,t})=\frac{1}{\sigma^2}\frac{\alpha}{1-e^{-\alpha t}},
$$
and $J_t(u^{*,t})$ increases to $J_T(u^*_\cdot)$ with $t\nearrow
T$. \qed

\section{\bf LDP in Discrete Time }\label{LDPdisc}

Let $m=\inf\{n\le N:u_n=0\}$ and $m=\infty$ if all $(u_n)_{n\le
N}$ are positive.
\begin{equation*}
I(y,x)=\sup_{t\in(-\infty,t_0)}[ty-x\frak{g}(t)].
\end{equation*}
\begin{theorem}\label{theo-A.1}
Assume \eqref{log}. For any $N\ge 1$, the family
$\big\{(x^K_n)_{n\le N}\}_{K\to\infty}$ obeys the LDP in
$\mathbb{R}^N_+$, supplied by the Euclidian metric $\varrho_N$,
with the speed $\frac{1}{K}$ and the rate function
$$
J_N(u_\cdot)= \left\{
\begin{array}{lll}
\sum\limits_{n=1}^{m-1} I(u_n,u_{n-1})-u_{m-1}\log(p_0), &
\substack{u_0=1\\u_n=0, n>m}
\\
\sum\limits_{n=1}^{N} I(u_n,u_{n-1}), & \substack{u_0=1\\u_n>0,
n\le N}.
\\
\infty, &\substack{\exists \ n:  u_n=0,u_{n+1}>0\\ \text{or} \
u_0\ne 1}
\\
\end{array}
\right.
$$
\end{theorem}

\begin{remark}
{\rm LDP for branching processes have been considered in the
literature, see, for example, \cite{Ath}, \cite{BigBin},
\cite{Ney}. However, they were concerned with the sequence
$\frac{X_n}{X_{n-1}}$, as $n\to\infty$, whereas here we consider
the LDP for $\frac{X_n}{X_0}$ processes indexed by the large
initial value. }
\end{remark}

\begin{remark}
{\rm The nonnegativity of $x^K_n$ provides some difficulty for
verification of LDP
 at the ``point of extinction''
where the test function becomes zero.  For set $\mathbb{S}$ of
  test functions that keep away from zero
the statement of the theorem  is implied by a result in Klebaner
and Zeitouni, \cite{KlZei} and other known results that can be
adapted to our setting (see, e.g. Kifer, \cite{Kif}, Puhalskii,
\cite{puh2}, Klebaner and Liptser, \cite{Klip},  etc.). But
$\{\tau\le N\}\not\in\mathbb{S}$ , and for the sake of
completeness and accuracy
 we give  the complete proof below, with a new   proof of the lower
bound in the local LDP.}
\end{remark}

\subsection{Proof of Theorem \ref{theo-A.1}}
We follow standard (necessary and sufficient) conditions for
proving the LDP by showing the exponential tightness:
\begin{equation*}
\lim_{C\to\infty}\varlimsup_{K\to\infty}\frac{1}{K}\log\mathsf{P}
\big(\varOmega\setminus\mathcal{K}_C\big)=-\infty
\end{equation*}
with compacts $\mathcal{K}_C=\{\max_{1\le n\le N}x_N\le C\}$,
$C\nearrow\infty$, and the local LDP:
\begin{equation*}
\lim_{\delta\to 0}\lim_{K\to\infty}\frac{1}{K}\log\mathsf{P}
\Big(\varrho_N(x^K_\cdot,u_\cdot)\le\delta\Big) =-J_N(u_.).
\end{equation*}

Notice that \eqref{log} implies the existence of a stochastic
exponential, with $t_n\le Kt_0$,
\begin{equation*}
\mathscr{E}^K_{(t_1,\ldots,t_N)}(x^K_1,\ldots,x^K_{N-1})
=\prod\limits_{n=1}^N
\mathsf{E}\Big(e^{t_nx^K_n}\big|\mathscr{F}_{n-1}\Big),
\end{equation*}
where $(\mathscr{F}_n)_{n\ge 0}$ is  the filtration, with
$\mathscr{F}_0=\{\varnothing,\varOmega\}$, generated by
$(x^K_n)_{n\ge 1}$.

Set
\begin{equation}\label{zzz}
\frak{z}_n=e^{\sum_{i\le n}t_\ell x^K_i -\log
\mathscr{E}^K_{(t_1,\ldots,t_n)}(x^K_1,\ldots,x^K_{n-1})}.
\end{equation}
The random process $(\frak{z}_n,\mathscr{F}_n)_{n\le N}$ is the
(positive) martingale,
\begin{equation}\label{=1}
\mathsf{E}\frak{z}_N=1.
\end{equation}

\subsubsection{Exponential tightness}
Since $ \max_{1\le n\le N}x^K_i\le\sum_{1\le n\le N}x^K_n$,
it is enough to show
\begin{equation*}
\lim_{C\to\infty}\varlimsup_{K\to\infty}\frac{1}{K}\log\mathsf{P}
\Big(\sum_{1\le i\le N}x^K_i\ge C\Big)=-\infty.
\end{equation*}

  Set $t^*=\argmax_{t\in(-\infty,t_0)}[t-\frak{g}(t)]$. Since
$\frak{g}(0)=0$, we have that $t^*\in(0,t_0)$ and $\frak{g}(t^*)<
t^*$. We choose $t_n\equiv t^*K(<Kt_0)$, and introduce
$\frak{A}=\big\{\sum_{1\le i\le n}x^K_i\ge C\big\}$. With chosen
$t_n$, we have $\mathsf{E}\frak{z}_N=1$ and, therefore,
$\mathsf{E}I_\frak{A}\frak{z}_N\le 1$. Taking into account this
inequality and \eqref{zzz}, write
\begin{equation*}
\begin{aligned}
1&\ge\mathsf{E}I_{\frak{A}}e^{\sum_{\{1\le n\le N\}}t^*x^K_n-\log
\mathscr{E}^K_{(t^*,\ldots,t^*)} (x^K_1,\ldots,x^K_{N-1})}
\\
&=\mathsf{E}I_{\frak{A}}e^{Kt^*\sum_{\{1\le n\le N\}}x^K_n-
K\frak{g}(t^*)\sum_{\{1\le n\le N\}}x^K_{n-1}}
\\
&\ge \mathsf{E}I_{\frak{A}}e^{K\sum_{\{1\le n\le
N\}}[t^*-\frak{g}(t^*)]x^K_n -K|\frak{g}(t^*)|}
\\
&\ge
\mathsf{E}I_{\frak{A}}e^{KC[t^*-\frak{g}(t^*)]}=e^{KC[t^*-\frak{g}(t^*)]
-K|\frak{g}(t^*)|} \mathsf{P}\big(\frak{A}\big).
\end{aligned}
\end{equation*}
 Therefore, $
\frac{1}{K}\log\mathsf{P}\big(\frak{A}\big)\le-\underbrace{[t^*-\frak{g}(t^*)]}_{>0}C
+|\frak{g}(t^*)| \xrightarrow[C\to\infty]{}-\infty. $ \qed

\subsubsection{Local LDP. Upper bound}
\label{sec-2.2.4}

We may restrict ourselves by the test function
$u_\cdot=\{\underbrace{u_1,\ldots,u_{N-1}}_{>0},\underbrace{u_N}_{=0}\}$
and show that
\begin{equation}\label{ubd}
\varlimsup_{\delta\to
0}\varlimsup_{K\to\infty}\frac{1}{K}\log\mathsf{P}
\Big(\varrho_N(x^K_\cdot,u_\cdot)\le\delta\Big) \le-J_N(u_\cdot).
\end{equation}
For the test function with all positive $u_n$'s and $u_0=1$ the
proof of \eqref{ubd} is similar. For test function with
$u_n=0,u_{n+1}>0$ or $u_0\ne 1$, \eqref{ubd} is obvious. For
others test functions the verification of \eqref{ubd} is reduced
to the above-mentioned ones.

Let now $\frak{A}=\big\{\rho_N(x^K_\cdot,u_\cdot)\le \delta
\big\}$. By \eqref{=1}, we have
\begin{equation}\label{aaa}
1\ge \mathsf{E}I_{\frak{A}}\frak{z}_N =
\mathsf{E}I_{\frak{A}}e^{\sum_{\{1\le n\le
N\}}[t_nx^K_i-Kx^K_{n-1}\frak{g}(\frac{t_n}{K})]}.
\end{equation}
Set $t^*_n=\argmax_{t\in(-\infty,t_0)}[tu_n-u_{n-1}\frak{g}(t)]$,
$n\le N-1$, and $t^*_N=-l \ (l>0)$, and take $t_n=Kt^*_n$, then we
derive from \eqref{aaa}
\begin{equation*}
\begin{aligned}
1&\ge \mathsf{E}I_{\frak{A}}e^{K\sum_{\{1\le n\le
N\}}[t^*_nu_n-u_{n-1}\frak{g}(t^*_n)] -K\sum_{1\le n\le
N-1}(t^*_n+|\frak{g}(t^*_n)|)\delta}
\\
& = \mathsf{E}I_{\frak{A}}e^{K[\sum_{\{1\le n\le
N-1\}}I(u_n,u_{n-1})-u_{N-1}\frak{g}(-l)] -K\sum_{1\le n\le
N-1}(|t^*_n|+|\frak{g}(t^*_n)|)\delta}
\\
& =
\mathsf{E}I_{\frak{A}}e^{K[J_{N-1}(u_\cdot)-u_{N-1}\frak{g}(-l)]
-K\sum_{1\le n\le N-1}(|t^*_n|+|\frak{g}(t^*_n)|)\delta}.
\end{aligned}
\end{equation*}
Hence, taking into account that
$\lim_{l\to\infty}\frak{g}(-l)=\log(p_0)$, we obtain
\begin{equation*}
\begin{aligned}
\frac{1}{K}\log\mathsf{P}\big(\frak{A}\big)&\le
-[J_{N-1}(u_\cdot)+u_{N-1}\frak{g}(-l)]+ \sum_{1\le i\le
N-1}(|t^*_i|+|\frak{g}(t^*_i)|)\delta
\\
&\xrightarrow[\delta\to
0]{}-[J_{N-1}(u_\cdot)+u_{N-1}\frak{g}(-l)]
\xrightarrow[l\to\infty]{}-J_{N}(u_\cdot).
\end{aligned}
\end{equation*}

\subsection{Local LDP. Lower bound}

Obviously for $u_\cdot$ with $J_N(u_\cdot)=\infty$, it is nothing
to verify. Further as in the upper bound verification, we may
restrict ourselves by the test function
$u_\cdot=\{\underbrace{u_1,\ldots,u_{N-1}}_{>0},\underbrace{u_N}_{=0}\}$
with $\mathsf{P}(\xi^1_1=0)=p_0>0$ and show that
\begin{equation*}
\varliminf_{\delta\to
0}\varliminf_{K\to\infty}\frac{1}{K}\log\mathsf{P}
\Big(\varrho_N(x^K_\cdot,u_\cdot)\le\delta\Big) \ge-J_N(u_\cdot).
\end{equation*}
Write
$$
\begin{aligned}
\{\varrho_N(x^k_\cdot,u_\cdot)\le\delta\}&=\big\{\varrho_{N-1}(x^k_\cdot,u_\cdot)
+x^K_N\le\delta \big\}
\\
& \supseteq \big\{\varrho_{N-1}(x^k_\cdot,u_\cdot)\le 0.5\delta,
\quad x^K_N\le 0.5\delta\big\}
\\
& \supseteq \big\{\varrho_{N-1}(x^k_\cdot,u_\cdot)\le
0.5\delta,\quad x^K_N=0\big\}
\\
& \supseteq \Big\{\varrho_{N-1}(x^k_\cdot,u_\cdot)\le 0.5\delta,
\quad \frac{1}{K}\sum_{j=1}^{Kx^K_{N-1} }\xi^j_N=0\Big\}
\\
& \supseteq \Big\{\varrho_{N-1}(x^k_\cdot,u_\cdot)\le
0.5\delta,\quad \frac{1}{K}\sum_{j=1}^{K(u_{N-1}
+\delta)}\xi^j_N=0\Big\}
\\
& = \Big\{\varrho_{N-1}(x^k_\cdot,u_\cdot)\le 0.5\delta,\quad
\sum_{j=1}^{K(u_{N-1} +\delta)}\xi^j_N=0\Big\}.
\end{aligned}
$$
The sets $\frak{A}_1=\big\{\varrho_{N-1}(x^k_\cdot,u_\cdot)\le
0.5\delta\big\}$ and $ \frak{A}_2=\big\{\sum_{j=1}^{K(u_{N-1}
+\delta)}\xi^j_N=0\big\} $ are independent, so that,
$$
\mathsf{P}\big(\varrho_N(x^k_\cdot,u_\cdot)\le \delta\big)\ge
\mathsf{P}\big(\varrho_{N-1}(x^k_\cdot,u_\cdot)\le 0.5\delta\big)
\mathsf{P}^{K(u_{N-1}+\delta)}\big(\xi^1_1=0\big).
$$
Consequently,
\begin{gather*}
\varliminf_{\delta\to 0}\varliminf_{K\to\infty} \frac{1}{K}\log
\mathsf{P}\big(\varrho_N(x^k_\cdot,u_\cdot)\le \delta\big)
\\
\ge \varliminf_{\delta\to 0}\varliminf_{K\to\infty}
\frac{1}{K}\log \mathsf{P}\big(\varrho_{N-1}(x^k_\cdot,u_\cdot)\le
0.5\delta\big)+ u_{N-1}\log\mathsf{P}\big(\xi^1_1=0).
\end{gather*}
If
\begin{equation}\label{lowbound}
\varliminf_{\delta\to
0}\varliminf_{K\to\infty}\frac{1}{K}\log\mathsf{P}
\Big(\varrho_{N-1}(x^{K}_\cdot,u_\cdot)\le\delta\Big)
\ge-J_{N-1}(u_\cdot),
\end{equation}
provided that $u_n>0$, $n\le N-1$, the   required lower bound
holds true.

Thus, it is left to verify the validity of \eqref{lowbound}.

Set $ \Lambda_{N-1}(x^K_\cdot)=\frak{z}_{N-1}, $ that is,
$$
\Lambda_{N-1}(x^K_\cdot)=e^{\sum_{n=1}^{
N-1}K\big[t^*_nx^K_n-x^K_{n-1}\frak{g}(t^*_n)\big]}, \quad
\mathsf{E}\Lambda_{N-1}(x^K_\cdot)=1.
$$
We introduce the probability measure $\mathsf{Q}^K_{N-1}$ with $
d\mathsf{Q}^K_{N-1}=\Lambda_{N-1}(x^K_\cdot)d\mathsf{P}. $ Since
$\Lambda_{n-1}(x^{K}_\cdot)>0$, $\mathsf{P}$-a.s., we also have
$
d\mathsf{P}=\Lambda^{-1}_{n-1}(x^K_\cdot)d\mathsf{Q}^K_{n-1}.
$

In particular, for $ \frak{A}=\big\{\varrho_{N-1}(
x^K_\cdot,u_\cdot)\le\delta\big\}, $
$$
\mathsf{P}(\frak{A})=\int_\frak{A}
\Lambda^{-1}_{N-1}(x^K_\cdot)d\mathsf{Q}^K_{N-1}.
$$
So, the following lower bound, on the set $\frak{A}$, is valid:
\begin{equation*}
\begin{aligned}
\Lambda^{-1}_{N-1}(x^K_\cdot)&\ge
e^{-KJ_{N-1}(u_\cdot)-K\delta\max_{n\le N-1}(|t^*_n|+
|\frak{g}(t^*_n)|)}
\\
&\ge e^{-KJ_{N-1}(u_\cdot)-K\delta\max_{n\le N-1}(|t^*_n|+
|\frak{g}(t^*_n)|)}
\end{aligned}
\end{equation*}
or, equivalently,
\begin{equation*}
\frac{1}{K}\log\mathsf{P}(\frak{A})\ge -J_{N-1}(u_\cdot)
-\delta\max_{n\le N-1}(|t^*_n|+|\frak{g}(t^*_n)|) +\frac{1}{K}\log
\mathsf{Q}^K_{N-1}(\frak{A}).
\end{equation*}
The latter inequality implies \eqref{lowbound} if
\begin{equation}\label{Q.a}
\lim_{K\to\infty}\frac{1}{K}\log \mathsf{Q}^K_{N-1}(\frak{A})=0.
\end{equation}
A  simple  condition, providing \eqref{Q.a}, is
$\lim_{K\to\infty}\mathsf{Q}^K_{N-1} (\frak{A})=1$ or,
equivalently,
\begin{equation}\label{oj}
\lim_{K\to\infty}\mathsf{Q}^K_{N-1}
\big(\varrho_{N-1}(x^K_\cdot,u_\cdot)>\delta\big)=0.
\end{equation}
We verify \eqref{oj} by showing\footnote{$\mathsf{E}^K_{N-1}$
denotes the expectation with respect to $\mathsf{Q}^K_{N-1}$}
\begin{equation}\label{posi}
\mathsf{E}^K_{N-1}\varrho^2_{N-1}(x^K\cdot,u_\cdot)=\frac{u_{N-1}}{K}\sum_{n=1}^{N-1}
\frac{u_{n-1}}{u^2_n}\frak{g}''(t^*_n).
\end{equation}
Notice that the positiveness of $(u_n)_{n\le N-1}$ provides a
boundedness for the right hand side of \eqref{posi}`     and, in
turn by Chebyshev's inequality, the validity of \eqref{oj}.

In order to establish \eqref{posi}, we apply the identity relative
to $t^*_n$:
\begin{equation}\label{4.10}
1=\mathsf{E}\Big(\frac{\Lambda_n(x^K_\cdot)}
{\Lambda_{n-1}(x^K_\cdot)}\Big|\mathscr{F}_{n-1}\Big)=\mathsf{E}e^{
  K\big[t^*_nx^K_n-x^K_{n-1}\frak{g}(t^*_n)\big]}.
\end{equation}
Differentiating twice  \eqref{4.10} in $t^{*}_n$, we find that
\begin{equation}\label{2.20a}
\begin{aligned}
0&=\mathsf{E}\Big([x^K_n-x^K_{n-1}
\frak{g}'(t^*_n)]\frac{\Lambda_i(x^K_\cdot)}
{\Lambda_{n-1}(x^K_\cdot)}\Big|\mathscr{F}_{n-1}\Big)
\\
0&=\mathsf{E}\Big(\big\{K[x^K_n-x^K_{n-1}
\frak{g}'(t^*_n)]^2-x^K_{n-1} \frak{g}''(t^*_n)
\big\}\frac{\Lambda_n(x^K_\cdot)} {\Lambda_{n-1}(x^K_\cdot)}\Big|
\mathscr{F}_{n-1}\Big).
\end{aligned}
\end{equation}
By the Bayes formula, e.g. \cite{LSII}, \cite{Kleb}: for any
integrable random variable $\alpha$,
$$
\mathsf{E}^K_{N-1}(\alpha|\mathscr{F}_{n-1})
=\mathsf{E}\Big(\alpha\frac{\Lambda_n(x^K_\cdot)}
{\Lambda_{n-1}(x^K_\cdot)}\Big| \mathscr{F}_{n-1}\Big).
$$
By taking $\alpha=x^K_n$ and $\alpha=[x^K_n-x^K_{n-1}
\frak{g}'(t^*_i)]^2$, we derive  with  the  help of  \eqref{2.20a}
that
\begin{gather}\label{4.12}
\mathsf{E}^K_{N-1}(x^K_n|\mathscr{F}_{n-1})=x^K_{n-1}\frak{g}'(t^*_n)
\\
\mathsf{E}^K_{N-1}\big([x^K_n-x^K_{n-1}
\frak{g}'(t^*_n)]^2|\mathscr{F}_{n-1}\big)=x^K_{n-1}\frac{\frak{g}''(t^*_n)}{K}
\label{4.13}
\end{gather}

Since $u_n,u_{n-1}$ are positive, we have
$\frak{g}'(t^*_n)=\dfrac{u_n}{u_{n-1}}$. Hence and by
\eqref{4.12}, we obtain that $
\mathsf{E}^K_{N-1}x^K_n=\frac{u_n}{u_{i-1}}\mathsf{E}^K_{N-1}x^K_{n-1}.
$ Consequently, iterating the above recursion and taking into
account $u_0=1$, we find that
$$
\mathsf{E}^K_{N-1}x^K_n=u_n.
$$
Further, with the help of  \eqref{4.13} we find a recursion
$$
\mathsf{E}^K_{N-1}(x^K_n)^2=\Big(\frac{u_n}{u_{n-1}}\Big)^2
\mathsf{E}^K_{N-1}(x^K_{n-1})^2+u_{n-1}\frac{\frak{g}''(t^*_n)}{K}.
$$
By using
$\mathsf{E}^K_{N-1}(x^K_n-u_n)^2=\mathsf{E}^K_{N-1}(x^K_n)^2-u^2_n$
and and $u^2_n=\big(\frac{u_n}{u_{n-1}}\big)^2u^2_{n-1}$, we
establish a recursion for
$\triangle_n=\mathsf{E}^K_{N-1}(x^K_n-u_n)^2$:
$$
\triangle_n=\Big(\frac{u_n}{u_{n-1}}\Big)^2\triangle_{n-1}
+u_{n-1}\frac{\frak{g}''(t^*_n)}{K}
$$
supplied by $\triangle_0=0$. Then, $\frac{\triangle_0}{u^2_0}=0$
and
$$
\frac{\triangle_n}{u^2_n}=\frac{\triangle_{n-1}}{u^2_{n-1}}
+\frac{u_{n-1}}{u^2_n}\frac{\frak{g}''(t^*_n)}{K}, \quad
\triangle_{N-1}=\frac{u_{N-1}}{K}\sum_{n=1}^{N-1}\frac{u_{n-1}}{u^2_n}
\frak{g}''(t^*_n).
$$
It is left to recall that $
\triangle_{N-1}=\mathsf{E}^K_{N-1}\varrho^2_{N-1}(x^K\cdot,u_\cdot).
$ \qed

\section{\bf LDP in Continuous Time}\label{LDPcont}

We introduce the filtration $(\mathscr{F}^B_t)_{t\ge 0}$ generated
by Brownian motion $B_t$, with the general conditions. All random
processes considered in this section are adapted to this
filtration. Henceforth,  by agreement,
$$
\frac{0}{0}=0.
$$

\begin{theorem}\label{theo-B.1}
For any $T>0$, the family $\big\{(x^K_t)_{t\le T}\}_{K\to\infty}$
obeys the LDP in $\mathbb{C}_{[0,T]}(\mathbb{R}_+)$, supplied by
the uniform metric $\varrho_T$, with the speed $\frac{1}{K}$ and
the rate function
$$
J_T(u_\cdot)= \left\{
\begin{array}{lll}
\frac{1}{2\sigma^2}\int_0^T\frac{(\dot{u}_t-\alpha u_t)^2}{u_t}dt,
& u_0=1, du_t=\dot{u}_tdt,
\\
\infty, & \text{otherwise}.
\end{array}
\right.
$$
\end{theorem}

\begin{remark}
{\rm Since $u_t\ge 0$, Freidlin-Wentzell's rate function,
\cite{FW}, $ \frac{1}{2\sigma^2}\int_0^T\frac{(\dot{u}_t-\alpha
u_t)^2}{u_t}dt $ is not compatible with $u_t=0$. Our branching
diffusion model is a very particular case of a model studied by
Puhalskii in \cite{PP3}. To apply the LDP analysis from \cite{PP3}
to the family $\big\{(x^K_t)_{t\le T}\}_{K\to\infty}$, one has to
``disentangle'' many details of the proof to make it compatible
with our case. Finally, in Donati-Martin et all, \cite{Yor}, the
LDP analysis deals with a rate function of the following type
$\int_0^T\frac{(\dot{u}_t-\rho)^2}{u_t}dt$ for $u_t\ge 0$ related
to a family of diffusion type processes without {\it extinction}.
A reader interested in details of the direct proof can find them
below. }
\end{remark}
\begin{proof}

It suffices to verify:

(i) $C$-exponential tightness (see \cite{LP}),
\begin{gather}
\label{B.1}
\lim_{C\to\infty}\varlimsup_{K\to\infty}\frac{1}{K}\log\mathsf{P}
\Big(\sup_{t\le T} x^K_t\ge C\Big)=-\infty,
\\
\lim_{\Delta\to 0}\varlimsup_{K\to\infty}\sup_{\gamma\le T}
\frac{1}{K}\log\mathsf{P}\Big(\sup_{t\le
\Delta}|x^K_{\gamma+t}-x^K_\gamma|\ge \eta\Big) =-\infty, \
\forall \ \eta>0, \label{B.2}
\end{gather}
where $\gamma$ is stopping time relative to $(\mathscr{F}^B_t)_{t\ge
0}$,

(ii) the Local LDP,
\begin{equation*}
\lim_{\delta\to
0}\lim_{K\to\infty}\frac{1}{K}\log\mathsf{P}\Big(\sup_{t\le T}
\big|x^K_t-u_t\big|\le\delta\Big)=-J_T(u_\cdot).
\end{equation*}

\medskip
{\bf (i)-Verification.} The It\^o equation \eqref{nd} is
equivalent to the integral equation
$
 x^K_t=e^{\alpha t}\big(1+\frac{1}{\sqrt{K}}\int_0^te^{-\alpha s}\sqrt{x^K_s}dB_s\big).
$
Hence,
\begin{equation}\label{OO1}
\sup_{t\le T}x^K_t\le 2e^{|\alpha|T}\Big(1\vee\frac{\sigma}{\sqrt{K}}\sup_{t\le T}
\int_0^te^{-\alpha s}\sqrt{x^K_s}dB_s\Big),
\end{equation}
so that, \eqref{B.1} holds true provided that
\begin{equation}\label{B.111}
\lim_{C\to\infty}\varlimsup_{K\to\infty}\frac{1}{K}\log\mathsf{P}\Big(\sup_{t\le T}
\int_0^te^{-\alpha s}\sqrt{x^K_s}dB_s\ge \sqrt{K}C\Big)=-\infty.
\end{equation}
In order to verify \eqref{B.111}, let us introduce a continuous martingale and its
variation process
\begin{equation*}
M_t=\frac{\sigma}{\sqrt{K}}\int_0^te^{-\alpha s}\sqrt{x^K_s}dB_s\quad\text{and}\quad
\langle M\rangle_t=\frac{\sigma^2}{K}\int_0^te^{-2\alpha s}x^K_sds
\end{equation*}
respectively and the stopping time $\tau_{C}=\inf\{t\le T:M_t\ge
C\}$, where $\inf\{\varnothing\}=\infty$ which enables us to claim
that \eqref{B.111} is valid if
\begin{equation}\label{tata}
\lim_{C\to\infty}\varlimsup_{K\to\infty}\frac{1}{K}\log\mathsf{P}\big(\tau_C\le T
\big)=-\infty.
\end{equation}
We proceed with  verification of \eqref{tata}. With $\lambda>0$,
set
$$
\frak{z}_t=e^{\lambda M_t-\frac{1}{2}\langle M\rangle_t}.
$$
It is well known that the process $(\frak{z}_t,\mathscr{F}^B_t)_{t\ge 0}$ is the
positive local martingale and so, the supermartingale too with
$\mathsf{E}\frak{z}_\theta\le 1$
for any stopping time  $\theta$ relative to $(\mathscr{F}^B_t)$.
By choosing $\theta=\tau_C$, we find that
$
1\ge \mathsf{E}I_{\{\theta\le T\}}\frak{z}_\theta.
$
Then, due to a lower bound on the set $\{\theta\le T\}$:
$
\log\frak{z}_\theta\ge \lambda C-\frac{\sigma^2\lambda^2}{2K}
\int_0^\theta e^{-2\alpha s}
x^K_sds
$
and \eqref{OO1} there exists positive $l$ such that
$$
\log\frak{z}_\theta\ge \lambda C-\frac{\sigma^2\lambda^2}{2K}
\int_0^\theta e^{-2\alpha s}\ge \lambda C-\frac{l\lambda^2}{2K}(1+C).
$$
Further, a choice of
$
\lambda=\frac{KC}{(1+C)l}
$
implies
$\frak{z}_\theta\ge e^{\frac{KC^2}{(1+C)l}}$.
Consequently,
$$
\frac{1}{K}\log\mathsf{P}\big(\tau_C\le T\big)\le -\frac{C^2}{(1+C)l}\xrightarrow
[C\to\infty]{}-\infty.
$$

\medskip
By \eqref{B.1}, the proof of \eqref{B.2} is reduced to the
verification of two conditions: for any $\eta,C>0$,
\begin{gather*}
\lim_{\Delta\to 0}\varlimsup_{K\to\infty}\sup_{\gamma\le T}
\frac{1}{K}\log\mathsf{P}\Big(\sup_{t\le
\Delta}\int_\gamma^{\gamma+t}x^K_sds\ge \eta, \ \sup_{s\le
T}x^K_s\le C\Big) =-\infty
\\
\lim_{\Delta\to 0}\varlimsup_{K\to\infty}\sup_{\gamma\le T}
\frac{1}{K}\log\mathsf{P}\Big(\frac{\sigma}{\sqrt{K}}\sup_{t\le
\Delta}\Big|\int_\gamma^{\gamma+t}\sqrt{x^K_s} dB_s\Big|\ge \eta,
\ \sup_{s\le T}x^K_s\le C\Big) =-\infty.
\end{gather*}
The first is obvious while the second is equivalent to
\begin{gather}\label{B.11}
\lim_{\Delta\to 0}\varlimsup_{K\to\infty}\sup_{\gamma\le T}
\frac{1}{K}\log\mathsf{P}\Big(\sup_{t\le \Delta}I_{T,C}
\big|M^K_{\gamma+t}-M^K_\gamma\big|\ge \eta\Big)=-\infty,
\end{gather}
where $I_{t,C}=I_{\{\sup_{s\le t}x^K_s\le C\}}$, $t\le T$.

Set $N^K_t=M^K_{\gamma+t}-M^K_\gamma$ and notice that
$(N^K_t,\mathscr{F}^B_{\gamma+t})_{t\ge 0}$ is a local martingale
with the variation process $\langle
N^K\rangle_t=\frac{\sigma^2}{K}\int_\gamma^{\gamma+t}x^K_sds $.

Further, the use of $ I_{T,C}N^K_t=I_{T,C}\int_0^tI_{s,C}dN^K_s $
simplifies \eqref{B.11}  to
\begin{gather}\label{B.12}
\lim_{\Delta\to 0}\varlimsup_{K\to\infty}\sup_{\gamma\le T}
\frac{1}{K}\log\mathsf{P}\Big(\sup_{t\le \Delta}
\Big|\int_0^tI_{s,C}dN^K_s\Big|\ge \eta\Big)=-\infty.
\end{gather}
The local martingale $N^{K,C}_t:=\int_0^tI_{s,C}dN^K_s$ possesses
the variation process
$$
\langle N^{K,C}\rangle_t=\int_0^tI_{s,C}d\langle N^K\rangle_s
=\frac{\sigma^2}{K}\int_0^tI_{s,C}x^K_sds,
$$
that is, d$\langle N^{K,C}\rangle_t\le \frac{\sigma^2C}{K}dt.$

Now, we are able to verify \eqref{B.12} with the help of
stochastic exponential technique. Let
$$
\frak{z}_t(\lambda)=e^{\lambda
N^{K,C}_t-\frac{\lambda^2}{2}\langle N^{K,C}\rangle_t}, \quad
\lambda\in \mathbb{R}.
$$
Since $\frak{z}_t(\lambda)$ is a continuous local martingale and
supermartingale too, for any stopping time $\theta$, $
\mathsf{E}\frak{z}_\theta(\lambda)\le 1. $ Let
$\theta=\inf\{t\le\Delta:N^{K,C}_t\ge \eta\}$. Taking into account
that $\{\theta\le \Delta\}=\{N^{K,C}_\theta\ge \eta\}$, write $
1\ge \mathsf{E}I_{\{\theta\le \Delta\}}\frak{z}_\theta(\lambda). $
The value $\frak{z}_\theta(\lambda)$ is evaluated below on the set
$\{\theta\le \Delta\}$ as follows: with $\lambda>0$ and
$
\langle N^{K,C}\rangle_\theta\le \frac{\sigma^2C}{K}\theta\le
\frac{\sigma^2C}{K}\Delta,
$
$$
\frak{z}_\theta(\lambda)\ge
e^{\lambda\eta-\frac{\lambda^2\sigma^2C}{2K}\Delta}.
$$
Therefore, $\log\mathsf{P}(\theta\le \Delta)\le-
\big[\lambda\eta-\frac{\lambda^2\sigma^2C}{2K}\Delta\big]$ and the
choice of $\lambda=\frac{K\eta}{\sigma^2C\Delta}$ provides
$$
\frac{1}{K}\log\mathsf{P}(\theta\le \Delta)\le
-\frac{\eta^2}{2\sigma^2C\Delta} \xrightarrow[\Delta\to
0]{}-\infty.
$$
It is clear that the same result remains valid for
$\theta=\inf\{t:-N^{K,C}_t\ge \eta\}$. Combining both, we
obtain \eqref{B.12}.

\medskip
{\bf (ii)-Verification. The upper bound.} For $u_0\ne 1$ or
$du_t\not\ll dt$,
  the proof is obvious.  For $u_0=1$ and
$du_t=\dot{u}_tdt$, the stochastic exponential technique is
applicable. With an absolutely continuous deterministic function
$\lambda(t)$ let us introduce a continuous martingale $M_t$ and
its predictable variation process $\langle M\rangle_t$:
$$
M_t=\frac{\sigma}{\sqrt{K}}\int_0^t\lambda(s)\sqrt{x^K_s} dB_s
\quad\text{and}\quad \langle
M\rangle_t=\frac{\sigma^2}{K}\int_0^t\lambda^2(s)x^K_s ds.
$$
It is well known that the stochastic exponential
$\frak{z}_t=e^{M_t-0.5\langle M\rangle_t}$ is a local martingale and
a supermartingale too with
$\mathsf{E}\frak{z}_T\le 1$. The use of this property implies
\begin{equation}\label{A.13}
1\ge \mathsf{E}I_{\{\sup_{t\le T}|x^K_t-u_t|\le\delta\}}\frak{z}_T.
\end{equation}
The next helpful step of the proof gives a \underline{deterministic} lower bound for
$\frak{z}_T$ on the set
$
\big\{\sup_{t\le T}|x^K_t-u_t|\le\delta\big\}=:\frak{A_\delta}.
$
By \eqref{KKK}, $M_t=\int_0^t\lambda(s)(dx^K_s-\alpha x^K_sds)$, so that,
\begin{equation*}
\begin{aligned}
\log\frak{z}_T&=\int_0^T\lambda(s)(dx^K_s-\alpha x^K_sds)-
\frac{\sigma^2}{2K}\int_0^t\lambda^2(s)x^K_s
ds
\\
&=\int_0^T\Big[\lambda(s)(\dot{u}_s-\alpha u_sds)-
\frac{\sigma^2}{2K}\lambda^2(s)u_s\Big]ds
\\
&\quad+\int_0^T\lambda(s)d(x^K_s-u_s)\quad [=\lambda_T(x^K_T-u_T)-\int_0^T(x^K_s-u_s)
\dot{\lambda}_tds]
\\
&\quad-\int_0^T\Big[\lambda(s)\alpha\big\{x^K_s-u_s\big\}+
\frac{\sigma^2}{2K}\lambda^2(s)\big\{x^K_s-u_s\big\}
\Big]ds.
\end{aligned}
\end{equation*}
Now, by taking $\lambda(s)=K\theta(s)$, we find a
lower bound of $\frak{z}_T$ on the set
$\frak{A}_\delta:=\{\sup_{t\le T}|x^K_t-u_t|\le \delta\}$,
\begin{equation*}
\begin{aligned}
\log\frak{z}_T&\ge K\int_0^T\Big[\theta(s)(\dot{u}_s-\alpha u_s)-
\frac{\sigma^2}{2}\theta^2(s)u_s\Big]ds
\\
&\quad
-\delta K\Big[|\theta_T|+\int_0^T\Big(|\dot{\theta}_s|+|\alpha\theta(s)|+
\frac{\sigma^2\theta^2(s)}{2}\Big)ds\Big].
\end{aligned}
\end{equation*}
This lower bound jointly with \eqref{A.13} implies
the following upper bound: for any absolutely continuous deterministic function
$\theta(s)$,
\begin{equation*}
\varlimsup_{\delta\to 0}\varlimsup_{K\to\infty}\frac{1}{K}\log
\mathsf{P}\big(\frak{A}_\delta\big)\le -
\int_0^T\Big[\theta(s)(\dot{u}_s-\alpha u_s)-
\frac{\sigma^2}{2}\theta^2(s)u_s\Big]ds.
\end{equation*}
Since $u_s$ is only nonnegative, it makes sense, for  computational convenience,
to use a corrected upper bound, with $\varepsilon>0$
\begin{equation}\label{B.14}
\varlimsup_{\delta\to 0}\varlimsup_{K\to\infty}\frac{1}{K}\log
\mathsf{P}\big(\frak{A}_\delta\big)\le -
\int_0^T\Big[\theta(s)(\dot{u}_s-\alpha u_s)-
\frac{\sigma^2}{2}\theta^2(s)(u_s+\varepsilon)\Big]ds.
\end{equation}
If $\dot{u}_t$ is absolutely continuous
function, a choice of $\theta(s)=\dfrac{\dot{u}_s-\alpha
u_s}{\sigma^2(u_s+\varepsilon)}$ provides
\begin{equation*}
\begin{aligned}
\varlimsup_{\delta\to 0}\varlimsup_{K\to\infty}\frac{1}{K}\log
\mathsf{P}\big(\frak{A}_\delta\big)&\le -\frac{1}{2\sigma^2}
\int_0^T\frac{(\dot{u}_s-\alpha u_s)^2}{u_s+\varepsilon}ds
\\
&\searrow -\frac{1}{2\sigma^2}
\int_0^T\frac{(\dot{u}_s-\alpha u_s)^2}{u_s}ds, \ \varepsilon\searrow 0.
\end{aligned}
\end{equation*}

In general case, one can choose a sequence $\theta_n(s)$, $n\ge 1$ of
absolutely continuous functions such that
\begin{gather*}
\lim_{n\to\infty}\Big[\theta_n(s)(\dot{u}_s-\alpha u_s)-
\frac{\sigma^2}{2}\theta^2_n(s)(u_s+\varepsilon)\Big]=
\sup_{\phi\in\mathbb{R}}\Big[\phi(\dot{u}_s-\alpha u_s)-
\frac{\sigma^2}{2}\phi^2(u_s+\varepsilon)\Big]
\\
=\frac{1}{2\sigma^2}\frac{(\dot{u}_s-\alpha u_s)^2}{u_s+\varepsilon}.
\end{gather*}
Hence, for sufficiently large $n$,
$
\big[\theta_n(s)(\dot{u}_s-\alpha u_sds)-
\frac{\sigma^2}{2}\theta^2_n(s)(u_s+\varepsilon)\big]\ge 0.
$
Then, due to \eqref{B.14} being valid with $\theta(s)$ replaced by $\theta_n(s)$, and
Fatou's theorem, we find that
\begin{gather*}
\varlimsup_{\delta\to 0}\varlimsup_{K\to\infty}\frac{1}{K}\log
\mathsf{P}\big(\frak{A}_\delta\big) \le -\varliminf_{n\to\infty}
\int_0^T\Big[\theta_n(s)(\dot{u}_s-\alpha u_s)-
\frac{\sigma^2}{2}\theta^2_n(s)(u_s+\varepsilon)\Big]ds
\\
\le -
\int_0^T\varliminf_{n\to\infty}\Big[\theta_n(s)(\dot{u}_s-\alpha u_s)-
\frac{\sigma^2}{2}\theta^2_n(s)(u_s+\varepsilon)\Big]ds
=-\frac{1}{2\sigma^2}\int_0^T\frac{(\dot{u}_s-\alpha u_s)^2}{u_s+\varepsilon}
\\
\searrow -\frac{1}{2\sigma^2}\int_0^T\frac{(\dot{u}_s-\alpha
u_s)^2}{u_s}, \ \varepsilon\searrow 0.
\end{gather*}

\medskip
{\bf (ii)-Verification.} The proof of
\begin{equation}\label{A.14}
\varliminf_{\delta\to 0}\varliminf_{K\to\infty}\frac{1}{K}\log
\mathsf{P}\Big(\sup_{t\le T}|x^K_s-u_s|\le \delta\Big)\ge
-\frac{1}{2\sigma^2}\int_0^T\frac{(\dot{u}_t-\alpha u_t)^2}{u_t}dt
\end{equation}
is done in three steps.

1. It suffices to analyse  the case $
\int_0^T\frac{(\dot{u}_s-\alpha u_s)^2}{u_s}ds<\infty $, which
enables us to consider only those  test functions that remain zero
after arriving at zero. In other words, we shall give the proof of
\eqref{A.14} for absolutely continuous $u_\cdot$ with $u_0=1$ and
$(u_t>0)_{t<T}$, $u_T\ge 0$.

\medskip
2. Set $ \tau_C=\inf\big\{t\le T:x^K_t\ge C\big\}, \quad
\text{where $\inf\{\varnothing\}=\infty $} $ and notice that if for any $C>0$
\begin{equation}\label{A.144}
\varliminf_{\delta\to 0}\varliminf_{K\to\infty}\frac{1}{K}\log
\mathsf{P}\Big(\sup_{t\le T\wedge\tau_C}|x^K_s-u_s|\le \delta\Big)\ge
-\frac{1}{2\sigma^2}\int_0^T\frac{(\dot{u}_t-\alpha u_t)^2}{u_t}dt,
\end{equation}
  then   \eqref{A.14} holds.
This can be seen as follows. Since
\begin{gather*}
\frak{A}_\delta\supseteq
\Big\{\sup_{t\le T\wedge\tau_C}|x^K_t-u_t|\le\delta\Big\}\cap\{\tau_C=\infty\}
\\
=\Big\{\sup_{t\le T\wedge\tau_C}|x^K_t-u_t|\le\delta\Big\}\setminus
\Big\{\sup_{t\le T\wedge\tau_C}|x^K_t-u_t|\le\delta\Big\}\cap\{\tau_C\le T\},
\end{gather*}
we have
$
\{\tau_C\le T\}\cup\frak{A}_\delta\supseteq
\big\{\sup_{t\le T\wedge\tau_C}|x^K_t-u_t|\le\delta\big\}
$, so that,
$$
2\Big[\mathsf{P}\big(\frak{A}_\delta\big)\vee\mathsf{P}\big(\tau_C\le T\big)\Big]\ge
\mathsf{P}\Big(\sup_{t\le T\wedge\tau_C}|x^K_s-u_s|\le \delta\Big).
$$
Hence, due to \eqref{A.144},
\begin{gather*}
\varliminf_{\delta\to 0}\varliminf_{K\to\infty}\frac{1}{K}
\log\mathsf{P}\big(\frak{A}_\delta\big)\bigvee\lim_{C\to\infty}
\varliminf_{K\to\infty}\frac{1}{K}
\log\mathsf{P}\big(\tau_C\le T\big)
\\
\ge-\frac{1}{2\sigma^2}\int_0^T\frac{(\dot{u}_t-\alpha u_t)^2}{u_t}dt
\end{gather*}
and it is left to recall that $\{\tau_C\le T\}=\big\{\sup_{t\le T}x^K_t\ge C\big\}$
and to refer to \eqref{B.1}.

\medskip
3. By 1., $\int_0^T\dot{u}^2_tdt<\infty$.
%
We proceed with the verification of \eqref{A.144}.
Define a continuous martingale $M_t$ and its variation process $\langle M
\rangle_t$: with $\varepsilon>0$,
$$
M_t=\int_0^{t\wedge\tau_C}\sqrt{K}\frac{\dot{u}_s-\alpha x^K_s}
{\sigma\sqrt{x^K_s+\varepsilon}}
dB_s\quad
\text{and}\quad
\langle M\rangle_t=
\int_0^{t\wedge\tau_C}K\frac{(\dot{u}_s-\alpha x^K_s)^2}{\sigma^2(x^K_s+\varepsilon)}ds.
$$
By definition of $\tau_C$, we have
$
\langle M\rangle_T\le \frac{2K}{\sigma^2\varepsilon}\int_0^T
\big(\dot{u}^2_t+\alpha^2C^2\big)ds<\infty,
$
so that the stochastic exponential $(\frak{z}_t,\mathscr{F}^B_t,\mathsf{P})_{t\le T}$
with
$
\frak{z}_t=e^{M_t-0.5\langle M\rangle_t}
$
is a uniformly integrable martingale, $\mathsf{E}\frak{z}_T=1$.
We use the latter property to define a new probability measure $\bar{\mathsf{P}}$ on
$(\Omega,\mathscr{F}^B_T)$ by letting
$
d\bar{\mathsf{P}}=\frak{z}_Td\mathsf{P}
$
and apply
\begin{equation*}
\mathsf{P}\Big(\sup_{t\le T\wedge\tau_C}|x^K_t-u_t|\le\delta\Big)=
\int\limits_{\{\sup_{t\le T\wedge\tau_C}|x^K_t-u_t|\le\delta\}}\frak{z}^{-1}_Td
\bar{\mathsf{P}}
\end{equation*}
for verification of \eqref{A.144}. This approach heavily uses a
semimartingale description of the processes
$(x^K_t,\mathscr{F}^B_t, \bar{\mathsf{P}})_{t\le T}$ and
$(\frak{z}^{-1}_t,\mathscr{F}^B_t, \bar{\mathsf{P}})_{t\le T}$. We
begin with the process $(B_t,\mathscr{F}^B_t,
\bar{\mathsf{P}})_{t\le T}$. The random processes
$(B_t,\mathscr{F}^B_t,\mathsf{P})_{t\le T}$ and
$(\frak{z}_t,\mathscr{F}^B_t,\mathsf{P})_{t\le T}$ are continuous
martingales and, in particular,
$$
d\frak{z}_t=I_{\{\tau_C\ge t\}}\frak{z}_t\sqrt{K}\frac{\dot{u}_t-\alpha x^K_t}
{\sigma\sqrt{x^K_t+\varepsilon}}dB_t.
$$
Hence, the co-variation process for $\frak{z}_t$, $B_t$ is defined
as:
$$
\langle \frak{z}, B\rangle_t=\int_0^{t\wedge\tau_C}
\frak{z}_s\sqrt{K}\frac{\dot{u}_s-\alpha x^K_s}
{\sigma\sqrt{x^K_s+\varepsilon}}ds.
$$
It is well known (see, e.g. Ch. 4, \S 5 in \cite{LSMar}) that
the random process
$
(\bar{B}_t,\mathscr{F}^B_t,\bar{\mathsf{P}})_{t\le T}
$
with
\begin{equation*}
\bar{B}_t=B_t-\int_0^t\frak{z}^{-1}_sd\langle \frak{z},B\rangle_s
=B_t-\int_0^{t\wedge\tau_C}\sqrt{K}\frac{\dot{u}_s-\alpha x^K_s}{\sigma(x^K_s+\varepsilon)}ds
\end{equation*}
is a Brownian motion. Consequently, we find that, $\bar{\mathsf{P}}$-a.s.,
\begin{equation}\label{bbb}
\begin{aligned}
x^K_t&=1+\int_0^tI_{\{\tau_C\ge s\}}\dot{u}_sds+\int_0^t\alpha
x^K_s \Big[1-I_{\{\tau_C\ge
s\}}\Big(\frac{x^K_s}{x^K_s+\varepsilon}\Big)^{0.5}\Big]ds
\\
&\quad
+\int_0^t\frac{\sigma}{\sqrt{K}}\sqrt{x^K_s}d\bar{B}_s
\\
\log{z}^{-1}_t&=-\int_0^{t\wedge\tau_C}
\sqrt{K}\frac{\dot{u}_s-\alpha x^K_s}{\sigma\sqrt{x^K_s+\varepsilon}}
d\bar{B}_s-\frac{1}{2}\int_0^{t\wedge\tau_C}K\frac{(\dot{u}_s-\alpha x^K_s)^2}
{\sigma^2(x^K_s+\varepsilon)}ds.
\end{aligned}
\end{equation}
Now, we evaluate from below the value $\frac{1}{K}\log\frak{z}^{-1}_T$ on the set
$\{\sup_{t\le T\wedge\tau_V}|x^K_t-u_t|\}$. Write
$$
\begin{aligned}
\frac{1}{K}\log\frak{z}^{-1}_T &\ge -\frac{1}{2\sigma^2}\int_0^{T}
\frac{(\dot{u}_s-\alpha u_s)^2}{u_s}ds
+h(C,\varepsilon,\delta)
\\
&\quad-\frac{1}{\sqrt{K}}\sup_{t\le T}
\Big|\int_0^{t\wedge\tau_C}
\frac{\dot{u}_s-\alpha x^K_s}{\sigma\sqrt{x^K_s+\varepsilon}}
d\bar{B}_s\Big|,
\end{aligned}
$$
where $h(C,\varepsilon,\delta)\xrightarrow[\delta\to 0]{}0$.
Therefore, \eqref{bbb} can be transformed into (here $\eta$ is a positive constant)
\begin{equation*}
\begin{aligned}
&\frac{1}{K}\log\mathsf{P}\Big(\sup_{t\le
T\wedge\tau_C}|x^K_t-u_t|\le\delta\Big) \ge
-\frac{1}{2\sigma^2}\int_0^{T} \frac{(\dot{u}_s-\alpha
u_s)^2}{u_s}ds +h(C,\varepsilon,\delta)
\\
&+\frac{1}{K}\log \bar{\mathsf{P}}\Bigg(\sup_{t\le T\wedge\tau_C}|x^K_t-u_t|\le\delta,
\frac{1}{\sqrt{K}}\sup_{t\le T}
\Big|\int_0^{t\wedge\tau_C}
\frac{\dot{u}_s-\alpha x^K_s}{\sigma\sqrt{x^K_s+\varepsilon}}
d\bar{B}_s\le \eta\Big|\Bigg).
\end{aligned}
\end{equation*}
This lower bound makes it possible to claim that \eqref{A.144}
holds true, provided that
\begin{gather}\label{5.16}
\lim_{K\to\infty}\bar{\mathsf{P}}
\Big(\sup_{t\le T}
\Big|\int_0^{t\wedge\tau_C}
\frac{\dot{u}_s-\alpha x^K_s}{\sigma\sqrt{x^K_s+\varepsilon}}
d\bar{B}_s\Big|>\sqrt{K}\eta\Big)=0
\\
\lim_{\varepsilon\to 0}\varlimsup_{K\to\infty}\bar{\mathsf{P}}
\Big(\sup_{t\le T\wedge\tau_C}|x^K_t-u_t|>\delta\Big)=0.
\label{5.17}
\end{gather}
Since $ I_{\{\tau_s\ge t\}}\frac{(\dot{u}_s-\alpha
x^K_s)^2}{\sigma^2(x^K_s+\varepsilon)}
\le\frac{(|\dot{u}_s|+C)^2}{\sigma^2\varepsilon} $ the Doob
inequality (here $\bar{\mathsf{E}}$ is the expectation relative to
$\bar{\mathsf{E}}$)
$$
\bar{\mathsf{E}}
\Big(\sup_{t\le T}
\Big|\int_0^{t\wedge\tau_C}
\frac{\dot{u}_s-\alpha x^K_s}{\sigma\sqrt{x^K_s+\varepsilon}}
d\bar{B}_s>\sqrt{K}\eta\Big|\Big)\le\frac{4}{K\eta^2}
\int_0^T\frac{(|\dot{u}_s|+C)^2}{\sigma^2\varepsilon}ds,
$$
jointly with 3., establish \eqref{5.16}.

Due to the first part of \eqref{bbb}, the proof of \eqref{5.17} is reduced to the
verification of
$$
\lim_{\varepsilon\to 0}\sup_{0\le x\le C}x
\Big[1-\Big(\frac{x}{x+\varepsilon}\Big)^{0.5}\Big]=0,
$$
which is obvious, and
$$
\lim_{K\to\infty}\bar{\mathsf{P}}
\Big(\sup_{t\le T}
\Big|\int_0^{t\wedge\tau_C}
\sqrt{x^K_s}d\bar{B}_s>\sqrt{K}\eta\Big|\Big)=0,
$$
which is similar to the proof of \eqref{5.16}.

\end{proof}

\end{document}